\documentclass[a4paper,10pt]{amsart}

\usepackage[colorlinks=true, linkcolor={black},citecolor={black}]{hyperref}

\usepackage{verbatim,exscale,amssymb,amsthm,mathrsfs,amsrefs}

\numberwithin{equation}{section}

\usepackage[all]{xy}

\xyoption{ps}
\SelectTips{cm}{}

\newdir{ >}{{}*!/-5pt/\dir{>}}

\newcommand{\even}{\textrm{ev}}
\newcommand{\supph}{\mathop{\mathrm{ssupp}}}

\newcommand{\eg}{\emph{e.g.}}
\newcommand{\ie}{\emph{i.e.}}
\newcommand{\cf}{\emph{cf.}}

\newcommand{\id}{\mathrm{id}}

\newcommand{\spec}{\mathop{\mathrm{Spec}}}
\newcommand{\specgr}{\mathop{\mathrm{Spec^h}}}
\newcommand{\Spec}{\mathop{\mathrm{Spec}}}
\newcommand{\spc}{\mathop{\mathrm{Spc}}}
\newcommand{\Spc}{\mathop{\mathrm{Spc}}}
\newcommand{\Proj}{\mathop{\mathrm{Proj}}}

\newcommand{\supp}{\mathop{\mathrm{supp}}}
\newcommand{\Supp}{\mathrm{Supp}}

\newcommand{\inthom}{\mathrm{\underline{Hom}}}

\newcommand{\Ann}{\mathrm{Ann}}

\newcommand{\op}{\mathrm{op}}

\newcommand{\ev}{\mathrm{ev}}

\newcommand{\colim}{\mathop{\mathrm{colim}}}

\newcommand{\unit}{\mathbf{1}} 
\newcommand{\obj}{\mathop{\mathrm{obj}}}

\newcommand{\proj}{\mathop{\mathrm{proj}}}

\newcommand{\Gr}{\mathrm{Gr}}

\newcommand{\Mod}{\mathrm{Mod}}
\newcommand{\GrMod}{\mathrm{GrMod}}

\newcommand{\Ch}{\mathrm{Ch}}

\newcommand{\D}{D} 

\newcommand{\Ab}{\mathrm{Ab}}

\newcommand{\Qcoh}{\mathop{\mathrm{Qcoh}}}

\newcommand{\FL}{R\text{-}\Gr\Mod_Z}

\newcommand{\Z}{\mathbb{Z}}

\DeclareMathOperator{\rank}{rank}

\theoremstyle{definition}

\newtheorem{defi}[equation]{Definition}

\newtheorem*{conv*}{Conventions}
\newtheorem*{ack}{Acknowledgements}
\newtheorem*{Related work}{Related work}

\theoremstyle{theorem}

\newtheorem{thm}[equation]{Theorem}
\newtheorem{lemma}[equation]{Lemma}
\newtheorem{thm-defi}[equation]{Theorem-Definition}
\newtheorem{prop}[equation]{Proposition}
\newtheorem{cor}[equation]{Corollary}

\newtheorem*{thm*}{Theorem}
\newtheorem*{lemma*}{Lemma}
\newtheorem*{cor*}{Corollary}
\newtheorem*{conj*}{Conjecture}
\newtheorem*{question*}{Question}

\theoremstyle{remark}

\newtheorem{notation}[equation]{Notation}
\newtheorem{remark}[equation]{Remark}
\newtheorem*{remark*}{Remark}

\newtheorem{example}[equation]{Example}
\newtheorem{examples}[equation]{Examples}

\begin{document}

\title{On the derived category of a graded commutative noetherian ring}

\author{Ivo Dell'Ambrogio}
\address{Universit\"at Bielefeld, Fakult\"at f\"ur Mathematik, BIREP Gruppe, Postfach 10\,01\,31, 33501 Bielefeld, Germany.}
\email{ambrogio@math.uni-bielefeld.de}

\author{Greg Stevenson}
\address{Universit\"at Bielefeld, Fakult\"at f\"ur Mathematik, BIREP Gruppe, Postfach 10\,01\,31, 33501 Bielefeld, Germany.}
\email{gstevens@math.uni-bielefeld.de}



 \date{}
 
 \maketitle

\begin{abstract}
For any graded commutative noetherian ring, where the grading group is abelian and where commutativity is allowed to hold in a quite general sense, we establish an inclusion-preserving bijection between, on the one hand, the twist-closed localizing subcategories of the derived category, and, on the other hand, subsets of the homogeneous spectrum of prime ideals of the ring. We provide an application to weighted projective schemes.
\end{abstract}

\tableofcontents

\section{Introduction}
\label{intro}

In his 1992 paper \cite{neemanChr}, Amnon Neeman has shown that for a noetherian commutative ring~$R$ ``one has a complete and very satisfactory description of the spectral theory of its derived category.'' 
Indeed, after providing a correct proof of Hopkins' classification of the thick subcategories of $\D^b(R\textrm-\proj)$ by way of specialization closed subsets of $\spec R$, he proceeds to show that even for the unbounded derived category $\D(R)$ one has an orderly classification, now in the form of a bijection
\begin{displaymath}\left\{ \begin{array}{c}
\text{subsets of}\; \spec R 
\end{array} \right\}
\xymatrix{ \ar[r]<1ex>^{\tau} \ar@{<-}[r]<-1ex>_{\sigma} &} \left\{
\begin{array}{c}
\text{localizing}\; \text{subcategories of} \; \D(R) \\
\end{array} \right\} .
\end{displaymath}
This restricts to a bijection between specialization closed subsets on one side and smashing subcategories on the other. From this one can prove the telescope conjecture for~$\D(R)$, which is the statement that every smashing subcategory of~$\D(R)$ (a localizing subcategory whose inclusion has a coproduct preserving right adjoint) is generated by the compact objects it contains. Another proof of the classification for $D^b(R\textrm-\proj)$ follows.

Neeman's bijection sends a subset $S\subseteq \spec R$ to the localizing subcategory $\tau(S)=\langle k(\mathfrak p) \mid \mathfrak p\in S\rangle $ generated by the residue fields at the primes in~$S$, and it sends a localizing subcategory $\mathcal L$ to the set $\sigma(\mathcal L)=\{\mathfrak p\in \spec R\mid \exists\, X\in \mathcal L \textrm{ s.t. } k(\mathfrak p) \otimes^\mathbb L X \neq 0\}$. In terms of the small support $\supph X= \{\mathfrak p\in \spec R \mid k(\mathfrak p) \otimes^\mathbb L X \neq 0 \}$ of a complex~$X$, it can be reformulated as follows:
\[
\tau(S)= \{X \in \D(R) \mid \supph X \subseteq S \}
\quad , \quad
\sigma(\mathcal L) = \bigcup_{X\in \mathcal L} \supph X
\,.
\]

More recently Dave Benson, Srikanth Iyengar and Henning Krause \cite{bik, bik3} have introduced a notion of support in the situation where one is given a compactly generated triangulated category~$\mathcal T$ together with an action by a $\Z$-graded commutative noetherian ring~$R$. 
In practice it is often the case that $\mathcal T$ is a rigidly-compactly generated tensor triangulated category, which we now assume (this means: $\mathcal T$ is also a tensor category with coproduct preserving exact tensor product~$\otimes$ and with compact unit object~$\unit$, and is such that its compact and rigid objects coincide). In this generality the support is given by
\[\supp{}_{\!R} \, X=\{\mathfrak p\in \specgr R \mid \mathit{\Gamma}_\mathfrak p(\unit) \otimes X \neq 0\}\] 
for each $X\in \mathcal T$, where $\specgr R$ is the homogeneous spectrum of~$R$, and where the $\mathit{\Gamma}_\mathfrak p(\unit)$ are suitable tensor idempotent objects provided by the theory.

A salient feature of the abstract theory is the identification of the following two hypotheses which together guarantee that the support $\supp_R$  provides a classification of localizing tensor ideals (here $\langle \mathcal F \rangle_{\otimes}$ denotes the localizing tensor ideal generated by a family of objects $\mathcal F\subseteq \mathcal T$): 
\begin{enumerate}
\item \emph{The (tensor) local-to-global principle}: $ \langle X\rangle_{\otimes}=\langle \mathit{\Gamma}_\mathfrak p(\unit) \otimes X \mid \mathfrak p\in \specgr R\rangle_\otimes$ for each object $X\in \mathcal T$.
\item \emph{Minimality}: the localizing tensor ideal $\langle \mathit{\Gamma}_\mathfrak p(\unit)\rangle_\otimes$ is either zero or minimal for each~$\mathfrak p\in \specgr R$.
\end{enumerate}
In \cite{bik2}, this theory is used to give a classification of the localizing tensor ideals of the stable module category of a finite group.
This spectacular success notwithstanding, the methods of \emph{loc.\,cit}.\ are for some purposes unsatisfactory; for instance, it is not clear how to free oneself from the affine tyranny of the graded ring. One would also like to capture classifications whose parameter space is, say, a noetherian scheme.

Some progress in this direction has been made quite recently in the PhD thesis of the second author (see \cite{StevensonActions}). In this work, the Benson-Iyengar-Krause theory is categorified as follows: $\mathcal T$ is a general compactly generated category, which is endowed with an action -- in a very natural sense -- by a rigidly-compactly generated category~$\mathcal R$; this includes or course the case where $\mathcal R=\mathcal T$ acts on itself via its tensor product. Now the relevant parameter space is $\spc \mathcal R^c$, the spectrum (in the sense of Paul Balmer~\cite{balmer_prime}) of the rigid-compact objects of~$\mathcal R$, and the objects $\mathit{\Gamma}_x(\unit)$ ($x\in \Spc \mathcal R^c$) are provided by the abstract theory of generalized Rickard idempotents of Paul Balmer and Giordano Favi~\cite{balmer-favi:rickard}.

In this new setting the local-to-global principle always holds, provided the category $\mathcal R$ has a Quillen model and the space $\Spc \mathcal R^c$ is noetherian. Moreover, by its very construction the whole theory is perfectly compatible with the extant powerful methods of tensor triangular geometry~\cite{balmer:icm}.
Once again, if one also has minimality then one gets a classification from the ensuing support theory in terms of some subsets of $\spc \mathcal{R}^c$. 

Using Neeman's classification one already understands how this theory recovers the support theory of Benson, Iyengar, and Krause when the derived category of a noetherian ring acts on a compactly generated triangulated category. As the theory of Benson, Iyengar, and Krause works in the generality of a $\Z$-graded noetherian ring it is thus natural to consider what happens when one allows the derived category of graded modules over such a ring to act. This requires a computation of the spectrum of the compact objects of such a derived category and the main goal of the current article is to perform this computation. Of course it is in any case natural to ask if Neeman's classification extends, in the obvious way, to graded modules over graded rings. The answer is yes; this opens the door to further exploring such categories, an undertaking which seems to have many potential applications.
In Section~\ref{sec:appl} we provide one first application to derived categories of weighted projective schemes.

\begin{center}$***$\end{center}

We allow $R$ to be graded by any abelian group~$G$, possibly with torsion, and we allow it to be commutative up to any reasonable sign rule, which covers both the strictly commutative case as well as the usual graded commutative one (see Def.\,\ref{defi:sign}). 
If $G$ is non trivial, $\D(R)$ is not generated by the tensor unit and therefore our geometric methods require us to restrict attention to those localizing subcategories which are tensor ideals; they are the same as those which are closed under twists by arbitrary elements $g$ of~$G$. 
Our classification is a bijection as follows:
\begin{displaymath}\left\{ \begin{array}{c}
\text{subsets of}\; \specgr R 
\end{array} \right\}
\xymatrix{ \ar[r]<1ex>^{\tau} \ar@{<-}[r]<-1ex>_{\sigma} &} \left\{
\begin{array}{c}
\text{twist-closed localizing}\; \text{subcategories of} \; \D(R) \\
\end{array} \right\} 
\end{displaymath}
with corollaries similar to Neeman's (see Theorem~\ref{thm:classification}).
As a special case we reprove Neeman's original classification for ordinary ungraded rings. That is not to say we give what could be considered a new proof; the ideas involved are essentially the same. However, our approach makes it abundantly clear which parts of the argument are completely general and belong to the realm of tensor triangulated categories, and which parts are instead specific to (graded) commutative noetherian rings. 

Here is a sketch of the proof. First, we notice that $\D(R)$ has a model (Prop.~\ref{prop:model}) and acts on itself by its symmetric tensor product~$\otimes^\mathbb L$.
Then we consider the \emph{small support} defined by the graded residue fields~$k(\mathfrak p)$
\[
\supph X = \{ \mathfrak p\in \specgr R \mid k(\mathfrak p)\otimes^\mathbb L X \neq 0 \}
\]
and we establish its properties, using only the most elementary results of the theory of injective objects in the category of graded $R$-modules. In particular, the small support detects the objects of~$\D(R)$ (Cor.\,\ref{cor:detection}), it is compatible with the tensor product (Lemma~\ref{lem_tensor_formula}) and it behaves nicely with respect to compact objects (Lemmas \ref{lem_closed} and~\ref{lem_exists}). By a quite general criterion, this is enough to establish a canonical homeomorphism
\[
\specgr R \cong \spc \D(R)^c
\]
(Theorem~\ref{thm:spectrum}). 
Since $R$ is graded noetherian, this space is noetherian.  Hence the abstract theory can be applied in its full power, and it remains to verify minimality; this follows easily from the identity $\langle \mathit{\Gamma}_\mathfrak p(\unit)\rangle_\otimes = \langle k(\mathfrak p)\rangle_{\otimes}$ (Prop.~\ref{prop:gamma_kappa}) and the ``field object'' property of the residue field $k(\mathfrak p)$ (Lemma \ref{lem_resobject}).

\begin{conv*}
All categories are $\Z$-categories and all functors are $\Z$-linear.
\end{conv*}

\begin{ack}
The authors are grateful to Estanislao Herscovich for a few helpful comments, as well as a healthy session of sign checking, during which his remarkable super-algebraic powers proved quite instructive.
\end{ack}

\section{Definitions and basic results}
\label{sec:def_basic}

Let $G$ denote our grading group, which will always be assumed to be abelian and whose operation will be written additively.
By a \emph{graded ring}~$R$ we always mean a unital and associative ring graded by~$G$; in other words, $R$ comes together with a decomposition
\[
R = \bigoplus_{g\in G} R_g
\]
such that the multiplication satisfies $R_g \cdot R_h \subseteq R_{g+h}$ for all $g,h\in G$, and thus also  $1\in R_0$.
A \emph{\textup(left\textup) graded module} over~$R$ is an $R$-module~$M$ together with a decomposition $M=\bigoplus_{g\in G}M_g$ such that $R_gM_h\subseteq M_{g+h}$. 
We denote by $R\textrm-\Gr\Mod$ the category of graded $R$-modules and
degree-zero homomorphisms, \ie, those $R$-linear maps $f:M\to N$ such that $f(M_g)\subseteq N_g$ for all $g\in G$. 
As customary we will write $\deg m = g$ to indicate that the degree of $m$ is~$g$, that is, that $m\in M_g$.

If $M$ is an $R$-module and $g$ and element of~$G$, we write $M(g)$ for~$M$ \emph{twisted by~$g$}, that is, $M$ endowed with the new $G$-grading with components $M(g)_h:= M_{h+g}$. 
We say that an $R$-module~$M$ is \emph{graded free} if~$M$ is a sum of twists of~$R$.

\begin{defi}
\label{defi:companion}
The \emph{companion category} of~$R$, denoted $\mathcal C_R$, is the small $\Z$-category whose set of objects is $\obj(\mathcal C_R):=\{\underline g \mid g\in G\}$, whose morphism groups are given by $\mathcal C_R(\underline g,\underline h) := R_{h-g}$, and with composition given by restricting the multiplication of~$R$ to the appropriate homogeneous components:
\begin{align*}
R_{\ell - h} \times R_{h-g} \longrightarrow R_{\ell - g} 
\quad, \quad (r,s)\mapsto rs
\,.
\end{align*}
\end{defi}

\begin{lemma}[\cite{graded}*{Proposition I.1.3}]
\label{lemma:functorial_picture}
There is an equivalence between $R\textrm-\Gr\Mod$ and the additive functor category $\Ab^{\mathcal C_R}$, given by the functor
\[
R\textrm-\Gr\Mod \to \Ab^{\mathcal C_R}
\quad ,\quad
M \mapsto (\underline g \mapsto M_g) 
\, .
\]
A quasi-inverse is provided by
\[
 \Ab^{\mathcal C_R} \to  R\textrm-\Gr\Mod
 \quad , \quad 
 F \mapsto \bigoplus_{g\in G} F(\underline g) 
 \,,
\]
where we endow the abelian group $\bigoplus_{g\in G} F(\underline g)$ with the grading on evident display and with the $R$-action induced by functoriality. 
Under this equivalence, the functor $\mathcal C_R(\underline g,-)$ corepresented by~$\underline g$ corresponds to the free left $R$-module $R(-g)$. 
\qed
\end{lemma}

In the following, we make the identification $R\textrm- \Gr\Mod = \Ab^{\mathcal C_R}$ whenever convenient. It follows from this description that the category $R\textrm-\Gr\Mod$ is Grothendieck abelian, and $\{R(g) \mid g\in G\}$ is a set of projective generators. This second statement follows immediately from the Yoneda lemma:

\begin{lemma} 
\label{lemma:yoneda}
There is a natural isomorphism of abelian groups 
\[
R\textrm- \Gr\Mod (R(-g) , M) \stackrel{\sim}{\longrightarrow} M_{g}  
\quad ,\quad 
f \mapsto f(1)
\]
for every $M\in R\textrm-\Gr\Mod$ and $g\in G$.
\qed
\end{lemma}

Note that, under the identification of Lemma~\ref{lemma:functorial_picture}, the Yoneda embedding becomes the fully faithful functor 
\[
R : (\mathcal C_R)^\op \longrightarrow R\textrm - \Gr\Mod
\]
with $R(\underline g)= R(-g)$ on objects and which sends the morphism $r\in \mathcal C_R(\underline g,\underline h)$ to right multiplication with~$r$, seen as an $R$-linear map $R(r) : R(-h)\to R(-g)$.

\begin{defi}
\label{defi:sign}
Let $\epsilon:G\times G\to \Z/2$ be a symmetric $\Z$-bilinear map. We say that the graded ring~$R$ is \emph{$\epsilon$-commutative} if  $r\cdot s = (-1)^{\epsilon(\deg r,\deg s)}s\cdot r $ holds for all homogeneous elements $r$ and~$s$ in~$R$.
\end{defi}

\begin{examples} $\phantom{M}$
\begin{enumerate}
\item If $\epsilon$ is identically zero, then $\epsilon$-commutative just means commutative.
\item For $G=\Z$ the integers and $\epsilon: \Z\times \Z\stackrel{\cdot}{\to} \Z\to \Z/2$ the multiplication map modulo two, we recover the familiar notion of a ($\Z$-)graded commutative ring. For instance, the graded endomorphism ring of the tensor unit object in any reasonable tensor triangulated category will be such a ring, by the Eckmann-Hilton argument (see \cite{suarez-alvarez}).
\item Hovey, Palmieri and Strickland~\cite{hps} entertain the notion of a multigraded unital algebraic stable homotopy category, in which one has a finite number of generating ``spheres'' $S^1, \ldots, S^d$, thus giving rise to an $\epsilon$-commutative $\Z^d$-graded endomorphism ring of the tensor unit, where the signing form $\epsilon:\Z^d\times \Z^d\to \Z/2$ is given by 
\[
\epsilon((n_1, \ldots,n_d),(n_1', \ldots,n_d')) = n_1n_1'+\ldots + n_dn_d' \quad \mathrm{mod} \;2 \, .
\]   
\item
A commutative superalgebra (or supercommutative algebra) is an algebra graded over $\Z/2$ which is $\epsilon$-commutative for the multiplication map $\epsilon:\Z/2\times \Z/2\stackrel{\cdot}{\to} \Z/2$. 
\end{enumerate}
\end{examples}

We will need to use localization for such $\epsilon$-commutative rings and to consider their homogeneous spectra; essentially these are given by the obvious constructions, but let us be (at least a little) explicit about what is meant and let us make a few comments on what happens at this level of generality. 
We begin, as in the usual $\Z$-graded case, by defining the even part of such a ring.

\begin{defi}
Let $R$ be an $\epsilon$-commutative $G$-graded ring, as in Definition~\ref{defi:sign}. We define its \emph{even part}, written $R^\even$, to be the commutative $G$-graded ring with components
\[
(R^\even)_g :=  
\left\{
\begin{array}{ll}
R_g & \textrm{ if } \epsilon(g, h)=0 \textrm{ for all }h\in G \\
0 & \textrm{ otherwise }  
\end{array}
\right.
\]
and with multiplication restricted from $R$. 
We say a homogeneous element is \emph{even} if it belongs to the even part and we say it is $\emph{odd}$ if it is not even.
\end{defi}

\begin{remark}
Note that the bilinearity and symmetry of~$\epsilon$ imply that $R^\even$ is indeed a well-defined unital subring of~$R$, which moreover is commutative.
Note also that with this definition odd elements may still belong to the center; \eg, if $R_0$ is a commutative ring and if we endow $R:=R_0[x]/(x^2)$ with the usual $\Z$-grading where $\deg x=1$, then the strictly commutative ring~$R$ is also $\epsilon$-commutative for the product $\epsilon:\Z\times \Z\to \Z/2$, for which $x$ is odd.
\end{remark}

Observe that when $R$ is $\epsilon$-commutative all homogeneous ideals are automatically two-sided. 
We denote by $\specgr R$  the collection of all homogeneous prime ideals of~$R$, and call it the  \emph{homogeneous spectrum} of~$R$. We will consider $\specgr R$ as a topological space with the Zariski topology. 
As usual, if we consider instead the homogeneous prime ideals of $R^\even$ we would get the same space, since the square of any homogeneous element is even. 
We say that $R$ is \emph{\textup(graded\textup) noetherian} if the ascending chain condition holds for homogeneous ideals of $R$. 

From this point forward all rings we consider are assumed to be noetherian. 

\begin{remark}
One may suspect that the noetherianity of $R$ should require the grading group~$G$ to be finitely generated. 
This is not the case: one can always artificially enlarge the grading group (extending $R$ by zero). A slightly less trivial example is if $R$ is a graded field (see the next section) in which case it is noetherian, independently of~$G$.
\end{remark}

We make the following easy observation about the homogeneous spectrum of such a ring.

\begin{lemma}\label{noethspec}
If $R$ is a noetherian $\epsilon$-commutative $G$-graded ring, then the spectrum of homogeneous prime ideals, $\specgr R$, is a noetherian topological space.
\qed
\end{lemma}

We can consider the graded localization of an $\epsilon$-commutative ring~$R$ at a multiplicative set~$S$ consisting of even (and therefore central) homogeneous elements. The construction of this localization is the obvious one and it enjoys the usual properties; in particular, it is again an $\epsilon$-commutative $G$-graded ring. 
Similarly, we can also localize any graded $R$-module at such a multiplicative subset. 
For a homogeneous prime ideal $\mathfrak{p}\subseteq R$ and a graded module~$M$, we denote by $M_\mathfrak{p}$ the localization of $R$ at the multiplicative set $S=R^\even \cap R^{\mathrm{h}} \cap (R \smallsetminus \mathfrak{p})$ of the even homogeneous elements of $R$ not in $\mathfrak{p}$. 
We observe that in this generality it is possible for odd elements to become invertible in such a localization. 

Next, we want to define a symmetric monoidal structure for graded $R$-modules, where $R$ is allowed to be any $\epsilon$-commutative $G$-graded ring. 
This can be done quite explicitly as follows.
Every left $R$-module $M={}_RM$ carries a canonical structure of right $R$-module, making it into an $R$-bimodule ${}_RM_{R}$, by setting
\begin{equation} \label{eq:right_action}
m \bullet r := (-1)^{\epsilon(\deg m, \deg r)} rm
\end{equation}
for all homogeneous $r\in R$ and $m\in M$. 
Every morphism of left $R$-modules is also a morphism of right modules for this action.
Then the tensor product $M\otimes_RN$ of $M$ with another left module~$N={}_RN$ is given in each component by the following quotient of abelian groups:
\[
(M\otimes_R N)_g := 
\frac{\bigoplus_{h} M_h \otimes_\Z N_{g-h} }{
\langle
 m \bullet r \otimes n -  m\otimes r n \\
\mid \; m\in M_{p}, r\in R_{h-p}, n \in N_{g-h} 
\rangle }
\]
(\cf\,\cite{graded}). 
The ring $R$ still acts on ${}_RM\otimes_R N$ on the left. There are evident natural associativity and right and left unit isomorphisms
\[
(L\otimes_R M)\otimes_R N\cong L\otimes_R (M\otimes_R N)
\quad, \quad
 M\otimes_R R \cong R
\quad, \quad
R \otimes_R M \cong R
\]
as well as a natural symmetry isomorphism:
\[
\tau_{M,N}:M \otimes_R N \cong N \otimes_R M
\quad, \quad 
\tau_{M,N}(m\otimes n):= (-1)^{\epsilon(\deg m, \deg n)} n\otimes m
\, .
\]

\begin{lemma}
\label{lemma:tensor}
The above constructions are well-defined and turn the category of graded left $R$-modules into a closed symmetric monoidal abelian category with tensor unit $R={}_RR$. 
\end{lemma}

\begin{proof}
All the verifications are straightforward and are therefore omitted. 
The existence of the internal Hom follows from the standard fact that every colimit preserving functor between Grothendieck categories, such as $M\otimes_R (-): R\textrm-\Gr\Mod \to R\textrm-\Gr\Mod $, has a right adjoint.  
\end{proof}

\begin{remark}
If one considers $R$ as a left $R$-module ${}_RR$, its canonical right action \eqref{eq:right_action} used for tensoring is just multiplication in~$R$ from the right, by $\epsilon$-commutativity. 
For the twisted left module $M={}_RR(g)$ however, beware that the element $m \bullet r$ in general is \emph{not} equal to the  product $m\cdot r$ computed in~$R$.
\end{remark}

\begin{lemma}
\label{lemma:twist_tensor}
There exist two natural isomorphisms 
\[
R(g) \otimes_R M \stackrel{\sim}{\longrightarrow} M(g)  
\quad \textrm{ and } \quad
M \otimes_R R(g) \stackrel{\sim}{\longrightarrow} M(g)
\]
of left $R$-modules for all $g\in G$ and all $M \in R\textrm-\Gr\Mod$.
\end{lemma}

\begin{proof}
The map $R(g)\otimes_R M\to M(g)$ given by $r \otimes m \mapsto (-1)^{\epsilon(g,\deg m)} rm$ is well-defined,  $R$-linear, and invertible with inverse $m\mapsto (-1)^{\epsilon(g,\deg m)} 1\otimes m$ (here $m\in M_{\deg m}= M(g)_{\deg m - g}$).
The second isomorphism is obtain by composing the first one with the switch isomorphism~$\tau_{R(g),M}$.
\end{proof}

If $M=R(h)$ in Lemma \ref{lemma:twist_tensor}, denote by 
\begin{displaymath}
\mu_{\underline{g}, \underline{h}} : \xymatrix{ R(-g) \otimes R(-h) \ar[r]^-{\sim} \ar@{=}[d] & R(-g-h) \ar@{=}[d] \\
R(\underline{g}) \otimes R(\underline{h}) \ar[r]^-{\sim} & R(\underline{g+h})
}
\end{displaymath}
the first isomorphism of left $R$-modules appearing in the lemma.

\begin{prop}
The companion category $\mathcal C_R$ of any $\epsilon$-commutative ring~$R$ carries a strict symmetric monoidal structure $\otimes$ with unit~$\underline 0$, given on objects by $\underline g \otimes \underline h := \underline{g+h}$ and on morphisms 
 $r\in \mathcal C_R(\underline g, \underline g')$ and 
 $s\in \mathcal C_R(\underline h, \underline h')$ by the formula 
\[
r\otimes s := (-1)^{\epsilon( g , h' -h )} rs
\; .
\]
With this tensor structure, the Yoneda embedding $R : (\mathcal C_R)^{\op}\to R\textrm-\Gr\Mod$ together with the identifications 
\[
\mu_{\underline g,\underline h}: R(\underline g)\otimes_R R(\underline h) \smash{\stackrel{\sim}{\longrightarrow}} R(\underline {g}\otimes \underline{h})
\quad \textrm{ and } \quad 
\mu_0 =\id :  R(\underline 0) \smash{\stackrel{=}{\longrightarrow}}  R
\]
becomes a strong symmetric monoidal functor $(R,\mu,\mu_0)$.
Moreover, the tensor category $\mathcal C_R$ is rigid.
\end{prop}

\begin{proof}
For the last assertion, note that the identity 
\[
\mathcal C_R(\underline g \otimes \underline h, \underline \ell) 
= R_{\ell -(g+h)}
= R_{(-h+\ell) -g}
= \mathcal C_R(\underline g, \underline {-h} \otimes \underline \ell)
\]
shows that each object $\underline h$ is rigid (\ie, strongly dualizable) with dual $\underline{-h}$.
All other verifications are straightforward and are therefore omitted.
\end{proof}

\begin{remark}
We stress that ``strict'' in the last proposition means that the associativity, left unit, right unit, and symmetry coherence isomorphisms are \emph{all} identity maps.
\end{remark}

\begin{remark}
It follows from the formal theory of Kan extensions --  or, in this context, Day convolution~\cite{day} -- that there exists, up to canonical isomorphism, a unique closed symmetric monoidal structure on the functor category $\Ab^{\mathcal C_R}$ such 
that the Yoneda embedding $(\mathcal C_R)^\op\to \Ab^{\mathcal C_R}$ is strong symmetric monoidal. 
By uniqueness we recover in this way the tensor product of Lemma~\ref{lemma:tensor}.
Indeed, this is how one can find the (rather quaint) formula for the tensor product in the companion category: given $r\in R_{g'-g}$ and $s\in R_{h'-h}$, one computes directly that the unique dotted map making the following square commute
\[
\xymatrix{
R(-g')\otimes_R R(-h') \ar[d]_{R(r)\otimes_R R(s)} \ar[r]^-{\mu_{\underline g',\underline h'}}_-\sim &
 R(-g'-h') \ar@{..>}[d] \\ 
R(-g)\otimes_R R(-h) \ar[r]^-{\mu_{\underline g,\underline h}}_-\sim & R(-g-h)
}
\]
is right multiplication by $(-1)^{\epsilon(g,h'-h)}rs$. Similarly one finds that the symmetry isomorphism $\tau_{R(g),R(h)}$ corresponds via $\mu$ to the identity map $R(g+h)\to R(h+g)$.
\end{remark}

Let $\Ch(R):= \Ch(R\textrm-\Gr\Mod)$ be the category of chain complexes of graded $R$-modules. 
It has a tensor product in the usual way, by setting
\[
(X \otimes_R Y)^n := \bigoplus_{p\in \Z} X^p \otimes_R Y^{n-p}
\quad \quad (n\in \Z)
\]
and by defining the differential with the Leibniz formula, for all complexes $X$ and~$Y$. 
The symmetric monoidal category $\Ch(R)$ is again closed, with the usual Hom complexes $\inthom_R(X,Y)$. 

\begin{prop}
\label{prop:model}
For every $\epsilon$-commutative  $G$-graded ring~$R$, the category $\Ch(R)$ has a (proper, cellular and combinatorial) Quillen model structure, where the weak equivalences are the quasi-isomorphisms and the fibrations are the degreewise surjections. Moreover, this model is compatible with the tensor product of complexes in the sense that it turns $\Ch(R)$ into a symmetric monoidal model category (see~\cite{hovey:model}).
\end{prop}

\begin{proof}
Of the various possibilities, we find it  most convenient to cite some results from~\cite{cisinski-deglise}. 
We recall that we have at hand a Grothendieck abelian category $\mathcal A := R\textrm-\Gr\Mod$ which is equipped with a closed symmetric monoidal structure. 
Moreover, it has a small set $\mathcal G := \{R(g) \mid g\in G\}$ of generators which contains the tensor unit $R=R(0)$ and which by Lemma~\ref{lemma:twist_tensor} is essentially closed under the tensor product.
It also follows immediately from Lemma~\ref{lemma:twist_tensor} that each $R(g)$ is flat, \ie, that the functor $R(g)\otimes_R( - ): \mathcal A\to \mathcal A$ is exact.

With this set $\mathcal G$ of generators, the \emph{$\mathcal G$-model structure} of \emph{loc.\,cit.\ }exists and has the properties listed in the proposition.
More precisely (and adopting the terminology of \emph{loc.\,cit.}), by~\cite{cisinski-deglise}*{Remark 1.15} it is always possible to choose a family $\mathcal H$ of complexes such that the pair $(\mathcal G,\mathcal H)$ forms a descent structure, so that by \cite{cisinski-deglise}*{Theorem 1.7} there exists a Quillen model structure on $\Ch(\mathcal A)$ -- which is independent of $\mathcal H$ other than for the choice of generating trivial cofibrations -- having quasi-isomorphisms for weak equivalences; the description of fibrations as the degreewise surjections (which will not be used in this article) follows from \cite{cisinski-deglise}*{Corollary 4.9} and the fact that every complex $X\in \Ch(\mathcal A)$ is $\mathcal G$-local,  that is, the canonical map
\[
K(\mathcal A) (\Sigma^n R(g), X) \to \D(\mathcal A) (\Sigma^n R(g), X)
\]
is bijective for all $n\in \Z$, where $ K(\mathcal A)$ denotes the homotopy category of complexes and $\Sigma$ the shift functor.
 
The facts that the generators~$\mathcal G$ are flat, include the tensor unit, and are essentially closed under tensoring, ensure that the model is compatible with the given symmetric monoidal structure, by \cite{cisinski-deglise}*{Proposition 2.8} and \cite{cisinski-deglise}*{Corollary 2.6}.
\end{proof}

\begin{remark}
Although the \emph{existence} of a model for $\D(R)$ will be required in Section~\ref{sec:spectrum_loc}, we will not have to actually work with it: the (probably) more familiar methods of homological algebra will amply suffice, see \eg\ \cite{keller}.
\end{remark}

It follows from Proposition \ref{prop:model} that, by deriving the tensor product and the internal Hom functors, the derived category of every $\epsilon$-commutative $G$-graded ring~$R$ inherits the structure of a closed tensor category $(\D(R), \otimes_R^\mathbb L , R, \mathbb R\inthom_R)$. 
Moreover, the tensor structure is compatible with the triangulation in the best way; we refer to \cite{hps}*{Appendix~A} for precise statements. 

If the group $\D(R)(\Sigma^n R(g), X)$ vanishes for all $n\in \Z$ and $g\in G$, then $X$ is acyclic. Hence 
\[
\{\Sigma^n R(g)\mid g\in G, n\in \Z\}
\]
is a set of compact generators for the triangulated category $D(R)$. 
It is not hard to see that the objects $\Sigma^n R(g)$ are also rigid, that is, that the canonical map
\[
\mathbb R\inthom_R(\Sigma^nR(g), R ) \otimes_R^\mathbb L X \to \mathbb R\inthom_R(\Sigma^nR(g), X)
\]
obtained by the tensor-Hom adjunction is an isomorphism for all~$X$. 
Hence $\D(R)$ is a \emph{rigidly-compactly generated tensor triangulated category}, as in \cite{hps} and \cite{balmer-favi:rickard}. 
In particular, the full subcategory $\D(R)^c\subseteq \D(R)$ of compact objects coincides with that of rigid objects.

\begin{notation}
Since no confusion should arise, we will simply write~$\otimes$ for the derived tensor product $\otimes_R^\mathbb L$ in~$\D(R)$. 
For any family  of objects $\mathcal F\subseteq \D(R)$ we will use the following notation:
\begin{align*}
\langle \mathcal F\rangle 
 & := \textrm{ the localizing subcategory of } \D(R) \textrm{ generated by } \mathcal F \, ,  \\
\langle \mathcal F\rangle_\otimes 
 & :=  \textrm{ the localizing tensor ideal  of } \D(R) \textrm{ generated by } \mathcal F \, .
\end{align*}
\end{notation}

The following observation will be used repeatedly.

\begin{lemma}
\label{lemma:tensor_ideals}
Let $R$ be any $\epsilon$-commutative $G$-graded ring and $\mathcal F\subseteq \D(R)$ any family of objects in the derived category. 
Then $\langle \mathcal F\rangle_\otimes$ coincides with the smallest localizing subcategory of $\D(R)$ containing $\mathcal F$ and closed under all the twist functors $(-)(g)$, $g\in G$, and also with the smallest localizing subcategory of $\D(R)$ containing the objects $\{X(g)\mid X\in \mathcal F, g\in G\}$.
\end{lemma}

\proof
The equality of the last two subcategories is obvious.
For the first one note that a localizing subcategory of~$\D(R)$ is a $\otimes$-ideal if and only if it is closed under tensoring with the generators $\Sigma^nR(g)$. 
It suffices therefore to show that there exist isomorphisms $R(g)\otimes X \cong X(g)$ for all complexes~$X\in \D(R)$, but this is an easy consequence of Lemma~\ref{lemma:twist_tensor}.
\qed

\section{Graded fields}
Fix an abelian group $G$ together with a $\Z/2\Z$-valued symmetric bilinear form~$\epsilon$. All rings considered henceforth are assumed to be $\epsilon$-commutative $G$-graded rings. Let us begin by recalling that a non-zero $\epsilon$-commutative $G$-graded ring $K$ is a \emph{graded field} if every non-zero homogenous element of $K$ is invertible. In particular, $K_0$ is an honest field, and the components~$M_g$ of every $K$-module~$M$ are $K_0$-vector spaces. We wish to show, in analogy with the ungraded case, that categories of modules over graded fields are rather structurally simple; this will provide us with a good theory of residue objects in the derived category, as in~\cite{neemanChr}.

We fix some graded field $K$ throughout the rest of the section.

\begin{defi}
Let $M$ be a graded $K$-module. We define the \emph{scaffold} of $M$ to be
\begin{displaymath}
s(M) := \{g\in G\; \vert \; M_g \neq 0\} \, .
\end{displaymath}
\end{defi}

\begin{lemma}
The scaffold, $s(K)$, of $K$ is a subgroup of $G$.
\end{lemma}
\begin{proof}
As $K$ is unital and $1\neq 0$ we must have $0\in s(K)$. If $g\in s(K)$ then there is a non-zero element in $K_g$ which, as $K$ is a graded field, must have an inverse in $K_{-g}$, so $-g\in s(K)$. Finally, suppose $g,g' \in s(K)$. Any non-zero element of $K_g$ gives, via multiplication, an isomorphism $K_{g'} \to K_{g+g'}$, so $g+g' \in s(K)$.
\end{proof}

\begin{lemma}
For any $g\in G$ we have $s(K(g)) = s(K) - g$.
\end{lemma}
\begin{proof}
Just note that
 $
s(K(g)) = \{h\in G \; \vert \; K_{g+h} \neq 0\} = s(K) - g.
 $
\end{proof}

\begin{lemma}
\label{lemma:free}
Every graded $K$-module $M$ is graded free.
\end{lemma}
\begin{proof}
If $h\in s(M)$ and $g\in s(K)$ then $g+h \in s(M)$ and $M_h \cong M_{g+h}$ because $K$ is a graded field. In particular the subgroup $s(K)$ of $G$ acts on $s(M)$ by translation. Let $\{h_i\}_{i\in I}$ be elements of $s(M)$ giving a decomposition of $s(M)$ into disjoint orbits $h_i + s(K)$. Then there is an isomorphism
\begin{displaymath}
\bigoplus_{i\in I} K(-h_i)^{m_i} \longrightarrow M \, ,
\end{displaymath}
where $m_i = \rank_{K_0} M_{h_i}$. Indeed, this is seen easily by choosing isomorphisms
\begin{displaymath}
K(-h_i)^{m_i}_{h_i} = K^{m_i}_0 \longrightarrow M_{h_i}
\end{displaymath}
and extending $K$-linearly.
\end{proof}

This gives us the next lemma, which is the graded analogue of \cite{bokstedt_neeman}*{Lemma 2.17}.

\begin{lemma}\label{lem_resobject}
Let $R$ be a $G$-graded ring and $R\to K$ a map (of graded rings) into a $G$-graded field~$K$. Then for all $X\in \D(R)$ the object $X\otimes K$ is a coproduct of suspensions and twists of $K$.
\end{lemma}
\begin{proof}
The functor $(-)\otimes K: \D(R)\to \D(R)$ factors through $\D(K)$ and so the result is immediate from Lemma~\ref{lemma:free}.
\end{proof}

\section{The small support}
Fix an abelian group $G$ and an $\epsilon$-commutative noetherian $G$-graded ring $R$. We now define a notion of support in terms of the graded residue fields of~$R$. 
We prove that this support satisfies all the desirable properties one would hope for. In this case the virtues of the support are not just their own reward: in the next section we see that a complete classification of the localizing $\otimes$-ideals of~$\D(R)$ follows in a very straightforward way from the results of this section and some abstract machinery. Let us begin by defining the objects which give rise to the small support.

\begin{defi}
Let $\mathfrak{p}\in \specgr R$ be a homogeneous prime ideal. We define the \emph{residue field} at $\mathfrak{p}$ in the usual way:
\begin{displaymath}
k(\mathfrak{p}) := R_\mathfrak{p}/\mathfrak{p}R_\mathfrak{p} = (R/\mathfrak{p})_{(0)} \,.
\end{displaymath}
\end{defi}
Happily it turns out that even in the $\epsilon$-commutative case this gives rise to graded fields.

\begin{lemma}
\label{lemma:residue_field_is_field}
Let $\mathfrak{p}$ be a homogeneous prime ideal of $R$. Then the residue field $k(\mathfrak{p})$ is a graded field. 
\end{lemma}
\begin{proof}
Let $r\in R_g$ be a homogeneous element of degree~$g$. 
Then $\deg r^2 = 2g$ and therefore $r^2$ is even. 
In particular, if $r\not\in\mathfrak p$ then $r^2\in (R\smallsetminus \mathfrak p)\cap R^\ev$ becomes inverted in~$k(\mathfrak p)$. 
But the inverse $r^{-2}$ is also even (of degree $-2g$). 
Therefore in $k(\mathfrak p)$ the element $r$ commutes with $r^{-2}$ and thus with $rr^{-2}$.
Hence $r$ is invertible with inverse~$rr^{-1}$.
\end{proof}

\begin{defi}
Let $X$ be an object of $\D(R)$. We define the \emph{small support} of $X$ to be the subset
\begin{displaymath}
\supph X := \{\mathfrak{p} \in \specgr R \; \vert \; k(\mathfrak{p}) \otimes X \neq 0 \textrm{ in }\D(R) \}
\end{displaymath}
of the homogeneous spectrum $\specgr R$ of $R$.
\end{defi}

\begin{remark}
We observe that there is no need to twist in this definition since, 
for every $g\in G$ and $X\in \D(R)$, we have 
$k(\mathfrak p)\otimes X \neq 0$ if and only if $ k(\mathfrak p)(g)\otimes X \neq 0$.
\end{remark}

Let us begin with those properties of the small support which are very obvious from the definition (so obvious in fact that we do not give a proof).

\begin{lemma}
\label{lem_evident}
The small support satisfies the following properties:
\begin{itemize}
\item[(i)] For every $X$ in $\D(R)$ we have $\supph X = \supph \Sigma X$.
\item[(ii)] For any set-indexed family $\{X_i\}_{i\in I}$ of objects of $\D(R)$ we have
\begin{displaymath}
\supph  \Big( \coprod_{i\in I} X_i \Big) = \bigcup_{i\in I} \supph X_i  \, .
\end{displaymath}
\item[(iii)] For any triangle $X\to Y \to Z \to \Sigma X$ in $\D(R)$ there is a containment
\begin{displaymath}
\supph Y \subseteq \supph X \cup \supph Z \, .
\end{displaymath}
\item[(iv)] $\supph R = \specgr R$.
\item[(v)] $\supph 0 = \varnothing$.
\end{itemize}
\end{lemma}

\begin{lemma}
\label{lem_tensor_formula}
The small support satisfies the tensor formula: for any $X$ and $Y$ in~$\D(R)$ we have
\begin{displaymath}
\supph (X\otimes Y) = \supph X \cap \supph Y \, .
\end{displaymath}
\end{lemma}
\begin{proof}
It is clear that $\supph(X\otimes Y)$ is contained in the intersection. To see the reverse inclusion just note that if $\mathfrak{p}$ is in the small support of both $X$ and $Y$ then
\begin{displaymath}
k(\mathfrak{p}) \otimes X \otimes Y 
\cong \Big(\coprod_i \Sigma^{m_i} k(\mathfrak{p})(g_i)^{\alpha_i} \Big) \otimes Y \neq 0 \, ,
\end{displaymath}
for some elements $g_i\in G$, integers $m_i$, and cardinals $\alpha_i$, where we use Lemma~\ref{lemma:residue_field_is_field}, Lemma~\ref{lem_resobject} and the fact that the tensor product commutes with coproducts.
\end{proof}

We next wish to check that the small support detects the vanishing of objects. This is, in some sense, the most technically unpleasant property to verify. However, most of the details are routine extensions of well known facts about $\Z$-graded rings to $G$-graded rings. 

The category $R$-$\GrMod$ is a locally noetherian Grothendieck abelian category. Thus it has enough injectives and every injective is a direct sum of indecomposable injectives. The general form of Matlis' theory (\cite{Matlis} and \cf\ \cite{Stenstrom}*{Chapter V.2}) shows that every indecomposable injective is the injective envelope, $E(R(g)/P)$, of $R(g)/P$ for some $g\in G$ and some irreducible submodule $P$ of $R(g)$. As twisting is an autoequivalence, it is easily seen that it is sufficient to consider only irreducible ideals of $R$. The argument of \cite{Matlis}*{Proposition 3.1} then extends in a straightforward way to show that every indecomposable injective is a twist of the envelope of $R/\mathfrak{p}$ where $\mathfrak{p}$ is a prime ideal, \ie, is of the form $E(R/\mathfrak{p})(g)$ (\cf\ \cite{BrunsHerzog}*{Theorem 3.6.3}). In particular, they are easily seen to be $\mathfrak{p}$-local and $\mathfrak{p}$-torsion in the graded sense.

\begin{prop}\label{prop_gen}
The objects $k(\mathfrak{p})(g)$, for $\mathfrak{p}\in \specgr R$ and $g\in G$, generate $\D(R)$:
\begin{displaymath}
D(R) 
= \langle\, k(\mathfrak{p})(g)\; \vert \; \mathfrak{p}\in \specgr R \; , \; g\in G \,\rangle = \langle\, k(\mathfrak{p})\; \vert \; \mathfrak{p}\in \specgr R \,\rangle_\otimes \,.
\end{displaymath}
\end{prop}
\begin{proof}
Let $X$ be a non-zero object of $\D(R)$ and pick $i\in \Z$ such that $H^i(X)\neq 0$. We may, without loss of generality, assume that $X$ is a complex of injectives by taking a K-injective resolution~(\cite{spaltenstein}).

Pick some non-zero homogeneous element of $H^i(X)$ and observe that it is represented by a morphism $f\colon \Sigma^{-i} R(g) \to X$ in $\D(R)$, which moreover may be assumed to correspond to a morphism of complexes, \ie, a map $R(g)\to X^i$. 
As $X^i$ is a direct sum of indecomposable injectives and $f$ is determined by the image of $1$, we may assume that the image of~$f$ is contained in a single indecomposable injective $E(R/\mathfrak{p})(g')\subseteq X^i$. Indeed $R(g) \to X^i$ factors through a finite direct sum of indecomposable injectives and if each of the restrictions of $f$ to these factors were null-homotopic, clearly $f$ would also be nullhomotopic. 
So we just replace~$f$, if necessary, by the restriction of~$f$ to a single indecomposable summand of $X^i$.

Since $E(R/\mathfrak{p})(g')$ is $\mathfrak{p}$-local, $f$ factors via $R_\mathfrak{p}(g)$. We can of course factor $R_\mathfrak{p}(g) \to E(R/\mathfrak{p})(g')$ through its image which is finitely generated over $R_\mathfrak{p}$ and $\mathfrak{p}$-torsion. Thus we get a factorization of $f$ via $(R_\mathfrak{p}/\mathfrak{p}^nR_\mathfrak{p})(g)$, for some integer~$n$. To summarize we have the following commutative diagram of factorizations of~$f$.
\begin{displaymath}
\xymatrix{
R(g) \ar[r] \ar[dr] \ar[d] & X^i \\
R_\mathfrak{p}(g) \ar[r] \ar[dr] & E(R/\mathfrak{p})(g') \ar[u] \\
& (R_\mathfrak{p}/\mathfrak{p}^nR_\mathfrak{p})(g) \ar[u]
}
\end{displaymath}
To complete the proof, just note that $(R_\mathfrak{p}/\mathfrak{p}^nR_\mathfrak{p})(g)$ is constructed from the $k(\mathfrak{p})(h)$, where $h\in G$, by taking finitely many extensions, so it certainly lies in the localizing subcategory generated by the $k(\mathfrak{p})(h)$. Hence, in~$\D(R)$, some $k(\mathfrak{p})(h)$ must also have a non-zero map to $X$.
\end{proof}

That the small support detects objects is an easy consequence of the proposition.

\begin{cor}
\label{cor:detection}
For every object $X$ of $\D(R)$ we have that $X\cong 0$ if and only if $\supph X = \varnothing$.
\end{cor}
\begin{proof}
One direction is clear. On the other hand, suppose $X\otimes k(\mathfrak{p})$ is zero for all $\mathfrak{p}\in \specgr R$. Then the kernel of the functor $X\otimes(-)$ is a localizing tensor ideal of $\D(R)$ containing all the residue fields. Hence it must be $\D(R)$ and it is immediate that $X\cong 0$.
\end{proof}

Before continuing, let us note the following important consequence of the last corollary.

\begin{prop}\label{prop_min}
For each $\mathfrak{p}\in \specgr R$, the localizing $\otimes$-ideal
\begin{displaymath}
\langle k(\mathfrak{p}) \rangle_\otimes = \langle k(\mathfrak{p})(g) \; \vert \; g\in G\rangle
\end{displaymath}
is minimal, \ie, it properly contains no non-zero localizing $\otimes$-ideal.
\end{prop}
\begin{proof}
Suppose $X \in \langle k(\mathfrak{p}) \rangle_\otimes$ is a non-zero object. 
Since $k(\mathfrak{p}) \otimes k(\mathfrak{q}) = 0$ whenever $\mathfrak{q}\neq \mathfrak{p}$, it must also hold that $X\otimes k(\mathfrak{q})$ is zero for all $\mathfrak q\in (\specgr R) \smallsetminus \{\mathfrak p\}$. 
By the last corollary we thus have that $X\otimes k(\mathfrak{p})$ is a non-zero object in $\langle X \rangle_\otimes$. 
It follows from Lemma \ref{lem_resobject} that the $\otimes$-ideal generated by~$X$ contains some twist of the object $k(\mathfrak{p})$, because localizing subcategories are thick. We conclude that $\langle X \rangle_\otimes = \langle k(\mathfrak{p}) \rangle_\otimes$.
\end{proof}

Next we wish to check that the small support of a compact object of $\D(R)$ is closed. 
For this we need the following lemma, whose ungraded analogue is well known.

\begin{lemma}  \label{lem_cpt_and_min} 
Let $R$ be a $\epsilon$-commutative noetherian $G$-graded ring.
\begin{enumerate}
\item[(i)]
An object is compact in $\D(R)$ precisely when it is isomorphic to a bounded complex of finitely generated projective graded modules. 

\item[(ii)] Let $(R,\mathfrak m, k)$ be graded local (e.g.\ $R_{\mathfrak p}$ for any homogeneous prime~$\mathfrak p$). 
Then in $\D(R)$ every right bounded complex of finitely generated projectives $C$ has a minimal graded free resolution $f:B\to C$; that is, $f$ is a quasi-isomorphism, the components~$B^i$ are finite graded free modules, and the differentials $d:B^i\to B^{i+1}$ satisfy $d(B^i)\subseteq \mathfrak m B^{i+1}$.
\end{enumerate}
\end{lemma}

\proof
(i) It is easily verified that bounded complexes of finitely generated projectives are compact; just use that, for such a complex $X$ and any other $Y\in D(R)$, one computes $\D(R)(X,Y)$ using homotopy classes of chain maps. 
To show the opposite inclusion note that, by the Thomason-Neeman localization theorem~\cite{neemanLoc}, $D(R)^c$ is the thick subcategory generated by the free modules $\{ R(g) \mid g\in G \}$.  
Hence it suffices to show that mapping cones and direct summands of bounded complexes of finite projectives are again of the same form; the first is clear, and the second follows (for instance) precisely as in \cite{buchweitz}*{Lemma 1.2.1}.

(ii)
The usual proof of the ungraded case, by induction, still works because Nakayama's lemma still holds for $G$-graded rings. 
\qed

\begin{lemma}\label{lem_closed}
Let $C$ be a compact object of~$\D(R)$. Then $\supph C$ is a closed subset of $\specgr R$.
\end{lemma}
\begin{proof}
Assume $C_\mathfrak p\neq0$ in $D(R_{\mathfrak p})$.
By Lemma~\ref{lem_cpt_and_min}, the complex $C_\mathfrak p$ has a minimal graded free resolution over~$R_\mathfrak{p}$. Tensoring with $k(\mathfrak{p})$ gives a complex which is non-zero in at least one degree and has zero differentials and so $C \otimes k(\mathfrak{p})$ is certainly non-zero. Conversely if $C_\mathfrak p=0$  then of course $C\otimes k(\mathfrak p)=C_\mathfrak p\otimes k(\mathfrak p)= 0$.
Hence $\supph C = V(\Ann_R H^*C) = \bigcup_i V(\Ann_R H^iC)$, which is closed since there are only finitely many non-vanishing cohomology groups. 
\end{proof}

Finally, we check that there are enough compact objects relative to the small support.

\begin{lemma}\label{lem_exists}
Let $V\subseteq \specgr R$ be a closed subset. Then there exists a compact object $C$ of $\D(R)$ such that $\supph C = V$.
\end{lemma}
\begin{proof}
By definition of the Zariski topology $V= V(I)$ for some homogeneous ideal $I\subseteq R$. 
Since $R$ is noetherian, we may write $I=(f_1, \ldots, f_n)$ for finitely many homogeneous elements  $f_i\in R_{g_i}$.
Let $C_i$ denote the mapping cone of~$f_i$, considered as a morphism $R(-g_i)\to R$. 
Each $C_i$ is a compact object, and therefore so is their tensor product $C:=C_1\otimes \cdots \otimes C_n$.
We claim that $\supph C = V$. 
Indeed, we have $\supph C= \supph C_1  \cap \cdots \cap \supph C_n$ by the tensor formula (Lemma~\ref{lem_tensor_formula}), so it suffices to show that $\supph C_i$ equals $V((f_i))$ for each~$i$.
By considering the triangle $R(-g_i)\to R\to C_i\to \Sigma R(-g_i)$, we see that $C_i\otimes k(\mathfrak p)\neq 0$ if and only if the morphism~$f_i \otimes k(\mathfrak p)$ is not invertible; that is, if and only if the element~$f_i$ belongs to the ideal~$\mathfrak p$. This proves the claim.
\end{proof}

\section{The spectrum and localizing tensor ideals}
\label{sec:spectrum_loc}

Let $\Spc D(R)^c$ denote the \emph{spectrum} of the tensor triangulated category of compact objects, in the sense of Balmer~\cite{balmer_prime}. 
We recall that this is a spectral topological space (defined for every essentially small tensor triangulated category) which comes together with a function $X\mapsto \supp X$ assigning a closed subset of $\Spc D(R)^c$ to every object $X\in D(R)^c$. The support function $\supp$ is compatible with the tensor triangular operations of $D(R)^c$, and it is the universal (finest) such.

\begin{thm}
\label{thm:spectrum}
For every $\epsilon$-commutative noetherian $G$-graded ring~$R$ there is a unique support preserving homeomorphism
\begin{displaymath}
\specgr R \to \spc \D(R)^c \,.
\end{displaymath}
In other words $(\specgr R, \supph)$ is a classifying support datum (\cite{balmer_prime}*{\S5}), meaning that there are inclusion preserving mutually inverse assignments
\begin{displaymath}
\left\{ \begin{array}{c}
\text{specialization closed} \\ \text{subsets of}\; \specgr R 
\end{array} \right\}
\xymatrix{ \ar[r]<1ex>^{\tau} \ar@{<-}[r]<-1ex>_{\sigma} &} \left\{
\begin{array}{c}
\text{thick}\; \otimes\text{-ideals } \\
\text{ of }  \; \D(R)^c\;
\end{array} \right\} 
\end{displaymath}
given, for a specialization closed subset $V$ of $\specgr R$ and a thick $\otimes$-ideal $\mathcal{J}$, by
\begin{displaymath}
\tau(V) = \{ X \in \D(R)^c \; \vert \; \supph X \subseteq V\}
\end{displaymath}
and
\begin{displaymath}
\sigma(\mathcal{J}) = \{\mathfrak{p}\in \specgr R \; \vert \; \exists X\in \mathcal J \textrm{ s.\,t. } \mathfrak p\in \supph X \} \, .
\end{displaymath}

\end{thm}
\begin{proof}
We wish to apply the recognition criterion \cite{kkGarticle}*{Theorem 3.1}. The category $\D(R)$ is rigidly-compactly generated and, by a serendipitous occurrence, we just happened to have proved  in the last section that $(\specgr R, \supph)$ satisfies all the necessary conditions to apply this criterion (by Corollary \ref{cor:detection} and Lemmas \ref{lem_evident}, \ref{lem_tensor_formula}, \ref{lem_closed} and \ref{lem_exists}).

The second part follows from the basic result of tensor triangular geometry, \cite{balmer_prime}*{Theorem~4.10}; in general, on the left hand side one would have to consider \emph{Thomason} subsets, but for the noetherian space $\specgr R$ these coincide with specialization closed subsets, and on the right hand side \emph{radical} thick $\otimes$-ideals, but since $D(R)^c$ is rigid these coincide with $\otimes$-tensor ideal, see \cite{balmer:supp}*{Proposition~2.4}.
\end{proof}

We now know that $\D(R)$ is a rigidly-compactly generated tensor triangulated category with a model and whose compacts have noetherian spectrum. Thus we can apply all of the machinery of~\cite{StevensonActions} to the problem of classifying the localizing $\otimes$-ideals of $\D(R)$. We shall mostly use this machinery, as well as the work of Balmer and Favi~\cite{balmer-favi:rickard}, as a black box; the following proposition spells out the little we need to know.

\begin{prop}\label{prop_ltg}
For each $x\in \spc \D(R)^c$ there exists a $\otimes$-idempotent object $\mathit{\Gamma}_xR$ of $\D(R)$ such that the assignment 
\begin{displaymath}
X\mapsto \supp X := \{x\in \spc \D(R)^c \; \vert \; \mathit{\Gamma}_xR\otimes X \neq 0\}
\quad (X\in \D(R))
\end{displaymath}
extends the Balmer support of compact objects and such that the following hold:
\begin{itemize}
\item[(i)] for $x\neq y$ we have $\mathit{\Gamma}_xR \otimes \mathit{\Gamma}_yR = 0$;
\item[(ii)] for every object $X$ of $R$ there is an equality of $\otimes$-ideals
\begin{displaymath}
\langle X \rangle_\otimes = \langle \mathit{\Gamma}_xR \otimes X \; \vert \; x \in \spc \D(R)^c \rangle_\otimes;
\end{displaymath}
\item[(iii)] an object $X$ of $\D(R)$ is zero if and only if $\supp X = \varnothing$.
\end{itemize}
\end{prop}
\begin{proof}
The construction and orthogonality of the idempotents is due to Balmer and Favi. Given the existence of a model and the fact that $\spc D(R)^c$ is noetherian the rest is a consequence of \cite{StevensonActions}*{Theorem~6.8}.
\end{proof}

\begin{lemma}
\label{lemma:unique_x}
For each $\mathfrak{p} \in \specgr R$ there is a unique $x\in \spc \D(R)^c$ such that $\mathit{\Gamma}_xR \otimes k(\mathfrak{p})$ is non-zero. In other words, $\supp k(\mathfrak{p}) = \{x\}$.
\end{lemma}
\begin{proof}
By part (iii) of the proposition we know there exists such an~$x$. Now suppose~$y$ is another point, distinct from ~$x$, such that $\mathit{\Gamma}_{y}R \otimes k(\mathfrak{p})$ is also non-zero. Then, using part~(i) of the proposition together with Lemma \ref{lem_resobject} we get
\begin{displaymath}
0 = \mathit{\Gamma}_{y}R \otimes \mathit{\Gamma}_{x}R \otimes k(\mathfrak{p}) 
\cong \mathit{\Gamma}_{y}R \otimes \Big(\coprod_i \Sigma^{n_i} k(\mathfrak{p})^{m_i}(g_i) \Big) \neq 0 \,,
\end{displaymath}
which is a contradiction.
\end{proof}

\begin{lemma}
If $\mathfrak{p} \neq \mathfrak{q}$ are two homogeneous prime ideals of $R$ then $k(\mathfrak{p})$ and $k(\mathfrak{q})$ have disjoint supports.
\end{lemma}
\begin{proof}
Suppose $k(\mathfrak{p})$ and $k(\mathfrak{q})$ both have support $\{x\}$. Then by the Half $\otimes$-Theorem (\cite{balmer-favi:rickard}*{7.22}) we have that, for any compact object $C$ of $D(R)$, 
\begin{align*}
\supp( k(\mathfrak p)\otimes C) =
\supp k(\mathfrak p) \cap \supp C =
\supp k(\mathfrak q) \cap \supp C =
\supp( k(\mathfrak q)\otimes C)
\, .
\end{align*}
Thus by Proposition \ref{prop_ltg} (iii) we see that $k(\mathfrak{p})\otimes C$ is zero if and only if $k(\mathfrak{q})\otimes C$ is zero. But this is clearly absurd as one can see, for example, from Lemma \ref{lem_exists}.
\end{proof}

\begin{prop}
\label{prop:gamma_kappa}
For every $\mathfrak p\in \specgr R$ there is an equality of localizing $\otimes$-ideals 
\begin{displaymath}
\langle \mathit{\Gamma}_xR \rangle_\otimes = \langle k(\mathfrak{p}) \rangle_\otimes 
\end{displaymath}
where $x$ is the unique point of $\spec \D(R)^c$ such that $\supp k(\mathfrak{p}) = \{x\}$. 
\end{prop}
\begin{proof}
The existence and uniqueness of $x$ is Lemma~\ref{lemma:unique_x}.
By \cite{StevensonActions}*{Proposition 5.5~(4)} we deduce from this that $\mathit{\Gamma}_xR \otimes k(\mathfrak{p}) \cong k(\mathfrak{p})$. 
On the other hand, it follows from the last lemma that $\mathit{\Gamma}_xR \otimes k(\mathfrak q)= 0$ for every $\mathfrak q\neq \mathfrak p$.
We know from Proposition~\ref{prop_gen} that the residue fields generate $\D(R)$ as a tensor ideal. 
So we have
\begin{align*}
\langle \mathit{\Gamma}_xR \rangle_\otimes &= \langle \mathit{\Gamma}_xR \rangle_\otimes \otimes \D(R) \\
&= \langle \mathit{\Gamma}_xR \rangle_\otimes \otimes \langle k(\mathfrak{q}) \; \vert \; \mathfrak{q}\in \specgr R \rangle_\otimes \\
&= \langle \mathit{\Gamma}_xR \otimes k(\mathfrak{q}) \; \vert \; \mathfrak{q}\in \specgr R \rangle_\otimes \\
&= \langle k(\mathfrak{p}) \rangle_\otimes
\end{align*}
where the third equality is an application of \cite{StevensonActions}*{Lemma 3.10}.
\end{proof}

\begin{cor}
For all $x\in \spc \D(R)^c$ the localizing $\otimes$-ideal $\langle \mathit{\Gamma}_xR \rangle_\otimes$ is minimal. Furthermore, the canonical homeomorphism of $\spc \D(R)^c$ with $\specgr R$ identifies $\supp X$ with $\supph X$ for all $X$ in $\D(R)$.
\end{cor}
\begin{proof}
The first statement is immediate from the proposition as the residue fields generate minimal $\otimes$-ideals by Proposition \ref{prop_min}. The second statement is a trivial consequence of the first.
\end{proof}

We can now easily deduce the classification theorem for localizing $\otimes$-ideals.

\begin{thm}
\label{thm:classification}
There are inclusion preserving mutually inverse bijections
\begin{displaymath}\left\{ \begin{array}{c}
\text{subsets of}\; \specgr R 
\end{array} \right\}
\xymatrix{ \ar[r]<1ex>^{\tau} \ar@{<-}[r]<-1ex>_{\sigma} &} \left\{
\begin{array}{c}
\text{localizing}\; \otimes\text{-ideals of} \; \D(R) \\
\end{array} \right\}, 
\end{displaymath}
and
\begin{displaymath}
\left\{ \begin{array}{c}
\text{specialization closed} \\ \text{subsets of}\; \specgr R 
\end{array} \right\}
\xymatrix{ \ar[r]<1ex>^{\tau} \ar@{<-}[r]<-1ex>_{\sigma} &} \left\{
\begin{array}{c}
\text{localizing}\; \otimes\text{-ideals of} \; \D(R)\; \\
\text{generated by objects of} \; \D(R)^c
\end{array} \right\} 
\end{displaymath}
where for a subset $W$ of $\specgr R$ and a localizing $\otimes$-ideal $\mathcal{L}$ we set
\begin{displaymath}
\tau(W) = \{ X \in \D(R) \; \vert \; \supph X \subseteq W\}
\end{displaymath}
and
\begin{displaymath}
\sigma(\mathcal{L}) = \{\mathfrak{p}\in \specgr R \; \vert \; \mathit{\Gamma}_\mathfrak{p}R \otimes \mathcal{L} \neq 0\}.
\end{displaymath}
\end{thm}
\begin{proof}
The map $\tau$ is a split monomorphism with left inverse $\sigma$ by \cite{StevensonActions}*{Proposition~6.3}. By the local-to-global principle (Proposition \ref{prop_ltg} (ii)) and its formal consequence \cite{StevensonActions}*{Lemma~6.2}, for every localizing $\otimes$-ideal $\mathcal{L}$ we have
\begin{align*}
\tau\sigma(\mathcal{L}) &= \tau(\{\mathfrak{p} \in \specgr R \; \vert \; \mathit{\Gamma}_\mathfrak{p}R\otimes \mathcal{L} \neq 0\})\\
&= \langle \mathit{\Gamma}_\mathfrak{p}R \; \vert \; \mathit{\Gamma}_\mathfrak{p}R \otimes \mathcal{L} \neq 0\rangle_\otimes
\, .
\end{align*}
To prove the first bijection note that, since the $\mathit{\Gamma}_\mathfrak{p}R$ generate minimal $\otimes$-ideals, we must have $\mathit{\Gamma}_\mathfrak{p}R \otimes \mathcal{L} = \langle \mathit{\Gamma}_\mathfrak{p}R \rangle_\otimes$ whenever this subcategory is non-zero. The result then follows from applying the local-to-global principle again.

The second pair of maps are well defined by Lemma \ref{lem_closed}, \cite{StevensonActions}*{Corollary 4.12}, and the good properties of the support. That they give a bijection is immediate from the first bijection.
\end{proof}

\begin{remark}
By Lemma~\ref{lemma:tensor_ideals}, one can reformulate the last theorem in the following way: there is an inclusion preserving bijection between subsets of $\specgr R$ and localizing subcategories of $D(R)$ closed under all twists $(-)(g)$.
\end{remark}

\begin{cor}
The category $\D(R)$ satisfies the relative telescope conjecture, \ie, the second bijection in Theorem~\ref{thm:classification} completely classifies those localizing $\otimes$-ideals whose inclusion admits a coproduct preserving right adjoint.
\end{cor}
\begin{proof}
This is an application of \cite{StevensonActions}*{Theorem 7.14}; one just needs to note that, by Lemma~\ref{lem_closed}, compact objects have closed supports (as we have identified our two notions of support), and that by Lemma \ref{lem_exists} any closed subset of $\specgr R$ can be realised as the support of a compact object.
\end{proof}

\section{An application to (weighted) projective schemes}
\label{sec:appl}
We now show how one easily obtains  from Theorem \ref{thm:classification} a classification of the localizing tensor ideals of the derived category of certain weighted projective schemes. In particular, if $R$ is a noetherian non-negatively $\Z$-graded ring this gives a direct method,  from an ``affine'' point of view, of classifying the tensor ideals in the derived category of $\Proj R$.

Let $R$ be a commutative noetherian $G$-graded ring, where $G$ is an abelian group, and as previously denote by $R$-$\Gr\Mod$ the category of graded $R$-modules. Let $Z$ be a closed subset of $\specgr R$ and denote by $U$ its open complement. We let $\FL$ denote the Serre subcategory of $R$-$\Gr\Mod$ consisting of those objects supported on $Z$ in the usual sense:
\begin{displaymath}
\FL = \{M \in R\text{-}\Gr\Mod \; \vert \; M_\mathfrak{p} = 0 \; \forall \; \mathfrak{p}\in U\}.
\end{displaymath}
We write 
\[
\Qcoh \mathbb{X} := R\textrm-\Gr\Mod / \FL
\]
for the abelian quotient of $R\textrm-\Gr\Mod$ by $\FL$. We think of $\Qcoh \mathbb{X}$ as the category of quasi-coherent sheaves on a ``weighted projective space~$\mathbb{X}$'' (precisely what this means is not important in the sequel, so let us not dwell on it). 
Observe that $\FL$ is the smallest Serre subcategory of $R$-$\Gr\Mod$ closed under filtered colimits and containing all twists of the residue fields of points in~$Z$. 

\begin{lemma}\label{lemma_htt}
The subcategory $\FL$ is the torsion class, $\mathcal{T}$, of the hereditary torsion theory on $R$-$\Gr\Mod$ cogenerated by
\begin{displaymath}
\{E(R/\mathfrak{p})(g)\; \vert\; \mathfrak{p}\in U, \; g\in G\} \,.
\end{displaymath}
\end{lemma}
\begin{proof}
Recall that for any $\mathfrak{p}\in \specgr R$ and $g\in G$ the indecomposable injective $E(k(\mathfrak{p}))(g)= E(R/\mathfrak p)(g)$ is $\mathfrak{p}$-torsion and $\mathfrak{p}$-local. Thus, by the universal property of localization, we must have $\FL \subseteq \mathcal{T}$.

On the other hand note that a finitely generated module lies in $\mathcal{T}$ if and only if its injective envelope is a (finite) direct sum of indecomposable injectives corresponding to points of~$Z$. As in the ungraded case one can easily check, using that filtered colimits of injectives in $R$-$\Gr\Mod$ are injective and localization preserves colimits, this extends to all objects of $\mathcal{T}$. Thus every object of $\mathcal{T}$ is a subobject of an object in $\FL$ (namely its injective envelope) and so $\mathcal{T} \subseteq \FL$ giving the claimed equality.
\end{proof}

It follows from the lemma that we have a diagram of abelian categories
\begin{equation}
\label{loc_ab}
\xymatrix{
\FL \ar[r]<0.5ex>^(0.4){} \ar@{<-}[r]<-0.5ex>_(0.4){} & R\text{-}\Gr\Mod \ar[r]<0.5ex>^(0.6){j^*} \ar@{<-}[r]<-0.5ex>_(0.6){j_*} & \Qcoh \mathbb{X}
}
\end{equation}
where the quotient~$j^*$ has right adjoint~$j_*$.

\begin{lemma}
\label{lemma:FL_tensor_ideal}
The subcategory $\FL$ is a tensor ideal, so that $\Qcoh \mathbb{X}$ inherits the tensor product of $R\textrm-\Gr\Mod$. 
\end{lemma}

\begin{proof}
This follows easily from the definition. Indeed, a module $M$ belongs to $\FL$ precisely when $\Supp_R M:= \{\mathfrak p\in \specgr R\mid M_\mathfrak p \neq 0\}$ is contained in~$Z$; now use that $\Supp_R(M\otimes N)\subseteq \Supp_R M \cap \Supp_R N$.  
\end{proof}

Let us now see what happens at the triangulated level. We denote by $D(R)^c_Z$ the thick subcategory of compact objects supported on~$Z$ (in the sense of Balmer). We let $\mathit{\Gamma}_ZD(R)$ be the localizing subcategory generated by $D(R)^c_Z$ and note that $\mathit{\Gamma}_ZD(R)$ is smashing as it is generated by compact objects of $D(R)$ (in fact it is precisely the subcategory $\tau(Z)$ as in Theorem \ref{thm:classification}). Let us be a little more explicit about what all this means. The subcategory $\mathit{\Gamma_Z}D(R)$ gives rise to a smashing localization sequence
\begin{equation}
\label{loc_seq}
\xymatrix{
\mathit{\Gamma}_ZD(R) \ar[r]<0.5ex>^(0.6){I_*} \ar@{<-}[r]<-0.5ex>_(0.6){I^!} & D(R) \ar[r]<0.5ex>^(0.5){J^*} \ar@{<-}[r]<-0.5ex>_(0.5){J_*} & L_Z D(R)
}
\end{equation}
\ie, all four functors are exact and coproduct preserving, $I_*$ and $J_*$ are fully faithful, $I^!$ is right adjoint to $I_*$, and $J_*$ is right adjoint to $J^*$. In particular there are associated coproduct preserving acyclization and localization functors given by $\mathit{\Gamma}_Z=I_*I^!$ and $\mathit{L}_Z=J_*J^*$ respectively. As in \cite{hps}*{Definition 3.3.2} this gives rise to Rickard idempotents $\mathit{\Gamma}_ZR$ and $L_ZR$ with the property that 
\begin{displaymath}
I_*I^! \cong \mathit{\Gamma}_ZR\otimes(-) \quad \text{and} \quad J_*J^* \cong L_ZR\otimes(-) \,;
\end{displaymath}
it follows that they are $\otimes$-orthogonal to each other by the usual properties of localization and acyclization functors. We observe that both $\mathit{\Gamma}_ZD(R)$ and $L_ZD(R)$ are tensor triangulated categories with units $\mathit{\Gamma}_ZR$ and $L_ZR$ respectively.

\begin{remark}
Such localization sequences are used to construct the idempotents $\Gamma_xR$ which appeared in Proposition \ref{prop_ltg}. More details can be found in \cite{balmer-favi:rickard}.
\end{remark}

\begin{lemma}
From the sequence of abelian categories \eqref{loc_ab} we obtain a localization sequence
\begin{displaymath}
\xymatrix{
D_{\FL}(R) \ar[r]<0.5ex>^(0.6){i_*} \ar@{<-}[r]<-0.5ex>_(0.6){i^!} & D(R) \ar[r]<0.5ex>^(0.5){j^*} \ar@{<-}[r]<-0.5ex>_(0.5){\mathbb Rj_*} & D(\Qcoh \mathbb{X}),
}
\end{displaymath}
where $D_{\FL}(R)$ denotes the full subcategory of $D(R)$ of complexes with cohomology in $\FL$.
\end{lemma}
\begin{proof}
First we show that the sequence in the statement is in fact a localization sequence. To prove this we need to check that $\mathbb{R}j_*$ is fully faithful and that the image of $i_*$ is the kernel of $j^*$.

We begin with the proof that $\mathbb{R}j_*$ is fully faithful. If $Y$ is an object of $D(\Qcoh \mathbb{X})$ then $j^*\mathbb{R}j_*Y$ is computed by taking a K-injective resolution $\tilde{Y}$ of $Y$ and applying $j^*j_*$. By \cite{Stenstrom}*{Chapter X, Proposition~1.4} an object of $\Qcoh \mathbb{X}$ is injective if and only if its image under $j_*$ is injective in $R$-$\Gr\Mod$. Thus $j_*\tilde{Y}$ is just this complex of injectives viewed in $D(R)$. In particular, $j^*\mathbb{R}j_*Y = j^*j_*\tilde{Y}$ is quasi-isomorphic to $Y$ via the natural map.

Let us now give the argument that $D_{\FL}(R)$ is the kernel of $j^*$. The functor $j^*$ is exact at the level of abelian categories and has kernel equal to $\FL$; thus $j^*$ commutes with taking cohomology and we see that its kernel consists precisely of those complexes whose cohomology modules lie in $\FL$.
\end{proof}

\begin{lemma}
The localization sequence of the last lemma agrees, up to monoidal equivalence, with~\eqref{loc_seq}.
\end{lemma}
\begin{proof}
First observe that $D_{\FL}(R)$ is a localizing $\otimes$-ideal of $D(R)$. It is clear that $D_{\FL}(R)$ is a localizing subcategory of $D(R)$ which is closed under twisting by all $g\in G$. Thus by Lemma \ref{lemma:tensor_ideals} it is a localizing $\otimes$-ideal.

So by Theorem \ref{thm:classification} the $\otimes$-ideal $D_{\FL}(R)$ must correspond to a subset of $\specgr R$ and this subset must contain $Z$ as $D_{\FL}(R)$ contains the residue field of each point in $Z$. It must in fact be $Z$ as if some $\mathfrak{q}\notin Z$ were in $\sigma (D_{\FL}(R))$ then, again by the classification, we would have $k(\mathfrak{q})$ in $D_{\FL}(R)$. But this is impossible since $k(\mathfrak{q})$ is not an object of $\FL$.
\end{proof}

\begin{cor}
There are inclusion preserving bijections
\begin{displaymath}\left\{ \begin{array}{c}
\text{subsets of}\; U 
\end{array} \right\}
\xymatrix{ \ar[r]<1ex>^{\tau} \ar@{<-}[r]<-1ex>_{\sigma} &} \left\{
\begin{array}{c}
\text{localizing}\; \otimes\text{-ideals of} \; \D(\Qcoh \mathbb{X}) \\
\end{array} \right\}, 
\end{displaymath}
and
\begin{displaymath}
\left\{ \begin{array}{c}
\text{specialization closed} \\ \text{subsets of}\; U 
\end{array} \right\}
\xymatrix{ \ar[r]<1ex>^{\tau} \ar@{<-}[r]<-1ex>_{\sigma} &} \left\{
\begin{array}{c}
\text{localizing}\; \otimes\text{-ideals of} \; \D(\Qcoh \mathbb{X})\; \\
\text{generated by objects of} \; \D(R)^c
\end{array} \right\}
\end{displaymath}
induced by the bijections of Theorem \ref{thm:classification}.
\end{cor}
\begin{proof}
By the last lemma it is sufficient to prove the result for $L_ZD(R)$. Note that any $\otimes$-ideal of $L_ZD(R)$ is also a $\otimes$-ideal in $D(R)$ as an object $X$ of $D(R)$ is in $L_ZD(R)$ if and only if it is isomorphic to $L_ZR \otimes X$. Thus the ideals in $L_ZD(R)$ are precisely those ideals of $D(R)$ contained in $L_ZD(R)$. Hence Theorem \ref{thm:classification} tells us that they are in bijection with subsets of $U$. The restricted bijection for those ideals generated by compact objects follows directly from the first bijection, the fact that the quotient functor to $L_ZD(R)$ sends compacts to compacts (see for instance \cite{NeeGrot}*{Theorem 5.1}), and \cite{StevensonActions}*{Lemma 7.10} which tells us that compacts of $L_ZD(R)$ have closed support in $U$.
\end{proof}

\begin{example}[Projective schemes]
Suppose $R$ is a non-negatively $\Z$-graded noetherian commutative ring which is generated by $R_1$ over~$R_0$. Then, letting $Z$ be the Zariski closure in $\specgr R$ of the irrelevant ideal~$R_{\geq1}$, $\Qcoh \mathbb{X}$ is equivalent to $\Qcoh (\Proj R)$, and by specializing the above result we see that the localizing tensor ideals of $D(\Qcoh (\Proj R))$ are in bijection with the subsets of $\Proj R$.
\end{example}

\begin{example}
Now suppose $R$ is a non-negatively $\Z$-graded finitely generated commutative $k$-algebra such that $R_0 = k$, where $k$ is some field. The grading on $R$ corresponds to an action of $k^*$ on $\spec R$ and the ideal $R_{\geq 1}$, generated by positively graded elements, corresponds to a closed fixed point $\mathbf{0}$ of the $k^*$ action. Letting $Z = \{\mathbf{0}\}$ there is an equivalence of categories
\begin{displaymath}
\Qcoh \mathbb{X} \simeq \Qcoh [(\Spec R \smallsetminus \mathbf{0})/k^*],
\end{displaymath}
where $[(\Spec R \smallsetminus \mathbf{0})/k^*]$ denotes the corresponding global quotient stack, as in \cite{Orlovgr}*{Proposition 28}. Thus we obtain a classification of localizing $\otimes$-ideals of the unbounded derived category of quasi-coherent sheaves on the quotient stack $[(\Spec R \smallsetminus \mathbf{0})/k^*]$ in terms of subsets of the punctured homogeneous spectrum. In particular, if $R$ is generated by $R_1$ over $R_0$ the quotient stack is just $\Proj R$ and we are in the situation of the previous example. See \cite{krishna} for related results.
\end{example}


\end{document}